\documentclass[a4paper,11pt]{amsart}

\usepackage[utf8]{inputenc}
\usepackage{mathrsfs}
\usepackage{amsfonts, amssymb, amsmath, amsthm}
\usepackage[active]{srcltx}
\usepackage{url}

\newcommand{\Q}{\mathbb{Q}}

\newcommand{\Z}{\mathbb{Z}}
\newcommand{\FF}{\mathbb{F}}

\usepackage[mathscr]{euscript}
\usepackage{latexsym}

\usepackage[all]{xy}

\usepackage{amssymb}
\usepackage{amsthm}
\usepackage{amsmath}
\usepackage{amsxtra}

.tfm
.tfm

\theoremstyle{plain}
\newtheorem{prop}{Proposition}[section]
\newtheorem{thm}[prop]{Theorem}
\newtheorem{coro}[prop]{Corollary}

\newtheorem*{property}{Property}

\theoremstyle{definition}

\theoremstyle{remark}

\newtheorem{remark}{Remark}

\newtheorem{example}[prop]{Example}

\numberwithin{table}{section}

\DeclareMathOperator{\Gal}{Gal}

\DeclareMathOperator{\sign}{sign}

\newcommand{\Om}{{\mathscr{O}}}

\def\ZZ{\mathbb Z}

\def\CC{\mathbb C}

\def\<#1>{{\left\langle{#1}\right\rangle}}

\def\Z{{\mathbb Z}}             
\def\Q{{\mathbb Q}}             
\def\QQ{{\mathbb Q}}             

\let\kro\dkro

\def\O{{R}}           
\def\Ol(#1){{\mathop{\O_l}(\id{#1})}}                 
\def\Or(#1){{\mathop{\O_r}(\id{#1})}}                 
\def\Orp(#1){{\mathop{\O_r}(\idp{#1})}}               
\def\Otern(#1){{\mathop{\O^0_r}(\id{a})}}    

\def\id#1{{\mathfrak{#1}}}      
\def\idp#1{{\id{#1}_\id{p}}}         

\def\T_#1(#2){{\mathop{\mathscr T}\nolimits_{#1}(\id{#2})}}
\def\TO(#1)_#2(#3){{\mathop{\mathscr T}\nolimits^{#1}_{#2}(\id{#3})}}

\def\A_#1(#2){{\mathop{\mathscr A}\nolimits_{#1}(\id{#2})}}
\def\Ax_#1(#2){{\mathop{\widetilde{\mathscr A}}\nolimits_{#1}(\id{#2})}}

\def\localconstant#1(#2){\varepsilon_{#1}({#2})}
\def\sign#1(#2){\epsilon_{#1}({#2})}

\DeclareMathOperator{\val}{val}
\DeclareMathOperator{\ab}{ab}
\DeclareMathOperator{\Ind}{Ind}
\DeclareMathOperator{\Tr}{Tr}
\DeclareMathOperator{\Norm}{Norm}
\DeclareMathOperator{\GL}{GL}
\DeclareMathOperator{\PGL}{PGL}
\DeclareMathOperator{\cond}{cond}
\begin{document}

\title{On the change of root numbers under twisting and applications}

\author{Ariel Pacetti}
\email{apacetti@dm.uba.ar}
\address{Departamento de Matem\'atica, Universidad de Buenos Aires,
         Pabell\'on I, Ciudad Universitaria. C.P:1428, Buenos Aires,
         Argentina}
\thanks{The first author was partially supported by PIP 2010-2012 GI and
UBACyT X867}


\begin{abstract}
The purpose of this article is to show how the root number of a modular
form changes by twisting in terms of the local Weil-Deligne
representation at each prime ideal. As an application, we show how one can
for each odd prime $p$, determine whether a modular form (or a Hilbert
modular form) with trivial nebentypus is Steinberg, Principal Series
or Supercuspidal at $p$ by analyzing the change of sign under a
suitable twist. We also explain the case $p=2$, where twisting is
not enough in general.
\end{abstract}

\maketitle

\section{Introduction}

The theory of local root number was developed in the well known work
\cite{AL}. Besides the definition and properties, in such paper the
authors also proved some cases of the variance of the local root
number under twisting. In particular, the results proven there imply
the well known result that the variance of the global root number of a
modular form of level $N$, by twisting by a quadratic field
corresponding to a character $\chi$ with conductor prime to $N$, is
given by $\chi(N)$. The case where the level and the conductor are not
prime to each other is more subtle. Some partial results were proven
in the same work, and some extensions with a similar perspective was
obtained in \cite{A-Li}.

The existence of local factors of representations was proved by
Deligne in \cite{Deligne}. Many authors used such description to
compute explicitly the local root numbers in terms of the local
representation of the Weil-Deligne group (as in \cite{KG}, or
\cite{Li}), but to our knowledge, although the variation of the local
factor under twisting follows essentially from the properties of the
local factor and is known to any expert in the area, it has not being
written down explicitly in general. As will be showed in the article,
it allows for example to easily compute for an elliptic curve the
local type (i.e. whether it is Steinberg, Principal Series or
Supercuspidal) at any odd prime (which of course can also be done by
looking at the reduction of the elliptic curve, as explained in
\cite{Roh}, Section 1). For $p=2$, a classification can be given in
terms of the sign variation, but this does not completely determine
the type (as is shown in Example \ref{ex:1} and Example
\ref{ex:2}). If one can have some other information, for example if
one is able to compute the space of modular forms appearing in the
quaternion algebra ramified at $2$, then this extra information is
enough to determine the type at $2$ as well.

It should be mentioned the recent article of \cite{LW}, were the
authors give a method to compute explicitly the local data of a
modular form (not just its type) for any nebentypus. They do not use
any twisting argument, so their method is different from ours. However,
they assume the modular form to be minimal between its twists (and it
is not clear how to compute the minimal twist without computing some
spaces of smaller levels with any nebentypus, see Example~\ref{ex:4}) and also to get the data
they need to compute the whole space of modular forms instead of just
the twist of a given one.

By the nature of the argument, the same results hold for Hilbert
modular form with trivial nebentypus as well, but the problem is that in
general the global characters to twist by, might not
exist. Nevertheless, this problem can be overcome by adding some
auxiliary prime to the twist, as is shown in the last section, so our
method works for Hilbert modular forms as well.

The current article started some years ago while studying a
characterization of the elliptic curves whose conductor is divisible
by $p$ but which do not show up in the quaternion algebra ramified at
$p$ and at $\infty$. After some numerical computations with Gonzalo
Tornar{\'\i}a we conjectured most of the formulas proven here, and
applied the formulas to find all elliptic curves in Cremona's tables
not appearing in any quaternion algebra (answering a question raised
by Professor Cremona to Gonzalo). The table can be found in
\url{http://mate.dm.uba.ar/~apacetti/}. For a work in progress with
Victor Rotger, we needed an exact formula giving the sign variation of
a modular form under twisting in terms of the local type, which forced
to reconsider the problem and by lack of references, write down a
proof of the conjectures.

I would like to thank many people who during the last year helped in a
way or another to this article: Gonzalo Tornar{\'\i}a, 
Tim Dokchitser, Lassina Demb\'el\'e, Victor Rotger, Luis
Dieulefait, Mladen Dimitrov and John Voight. Lassina and John provided
the examples of the last section, and taught me some useful things
about computing with Hilbert modular forms.
Also I would like to thank the people
from the CRM in Barcelona for the pleasant time I spent there and the
nice working environment.


{\bf Notation:} Given an odd prime number $p$, by $\val_p(N)$ we denote the
$p$-adic valuation of the integer $N$. By $p^\star$ we denote the
element $\kro{-1}{p}p$, whose square root generates the quadratic extension of
$\QQ$ ramified only at $p$. For a positive integer $r$, we denote by
$U_r$ the usual filtration in the $p$-adic units, given by
\[
U_r = \{x \in \Z_p^\times \, : \, x \equiv 1 \pmod{p^r}\}.
\]
If $E_p$ is a finite extension of $\QQ_p$, $\Tr$ and $\Norm$ denote
the trace and norm from $E_p$ to $\QQ_p$ respectively.

\section{Some useful well known results}

For $K$ a local field, let $W(K)$ denote the Weil group of $K$; that
is the preimage under the map
\[
\Gal(\overline{K}/K) \mapsto \Gal(\overline{k}/k),
\]
of the integer powers of Frobenius, where $k$ denotes the residue
field of $K$. Local class field theory gives an isomorphism between
$W(K)^{\ab}$ and $K^\times$. Furthermore, if $L/K$ is finite, the
following diagrams are commutative:

\begin{equation}
\xymatrix{
W(L) \ar[r] \ar@{^{(}->}[d]& W(L)^{\ab} \ar[d] \ar[r]^{\simeq}
&L^\times \ar[d]^{N_{L/K}}\\
W(K) \ar[r] & W(K)^{\ab} \ar[r]^{\simeq} &K^\times,}
\label{eq:commutative1}
\end{equation}
and
\begin{equation}
\xymatrix{
W(K)^{\ab} \ar[d]^{t} \ar[r]^{\simeq}& K^\times \ar@{^{(}->}[d]\\
W(L)^{\ab} \ar[r]^{\simeq}& L^\times,
}
\label{eq:commutative2}
\end{equation}
where $t$ denotes the transfer map.

Let $f \in S_k(\Gamma_0(N)$ be a non-zero weight $k$ and
level $N$ newform, i.e. a new modular form which is an eigenvalue for all
Hecke operators, 
and let $\rho_p(f)$ be the local representation of the Weil-Deligne
group $W'(\QQ_p)$ associated to $f$ at the prime $p$. Although it is
not so easy to give a description of the Weil-Deligne group (see
\cite{Tate}), it is relatively easy to describe its representations. A
complex $2$-dimensional representation of $W'(\QQ_p)$ is a pair
$(\rho,N)$, where:
\begin{enumerate}
\item $\rho$ is a representation $\rho:W(\QQ_p) \mapsto \GL_2(\CC)$,
\item $N$ is a nilpotent endomorphism of $\CC^2$ such that
\[
wNw^{-1} = \omega_1(w)n, \text{ for all } w \in W(\QQ_p),
\]
where $\omega_1$ is the unramified quasi-character giving the action of
$W(\QQ_p)$ on the roots of unity (and corresponds to the norm
quasi-character $||\, . \,||_p$ under local class field theory).
\end{enumerate}

Although representations of the Weil-Deligne group for any vector
space $V$ are defined in a similar way, the two dimensional complex
case is enough for our purposes.

\medskip



The correspondence between local components of automorphic
representation for $\Gamma_0(N)$ and the local representations of the Weil-Deligne
group, is given as follows (using the normalization given by Carayol
in \cite{carayol}):
\begin{enumerate}
\item {\bf Principal Series} (reducible case): the endomorphism $N = 0$ and 
\[
\rho_p(f) = \chi \oplus \chi^{-1}\omega_1^{1-k},
\]
for some quasi-character $\chi:W(\QQ_p)^{\ab} \mapsto \CC^\times$.
\item{\bf Steinberg or Special Representation} (indecomposable but
  reducible as $W(\QQ_p)$-representation): The endomorphism $N$ is
  given by the matrix $\left(\begin{smallmatrix}0 & 1\\ 0 &
    0 \end{smallmatrix} \right)$ and the representation $\rho_p(f)$ is
  given by
\[
\rho_p(f)(w) = \left( \begin{array}{cc} \chi(w) \omega_1(w) & 0\\ 0
  & \chi(w)\end{array} \right),
\]
for some quasi-character $\chi:W(\QQ_p) \mapsto \CC^\times$ with $\chi^2|_{\ZZ_p^\times}=1$.
\item{\bf Supercuspidal Representation I} (irreducible case, but
  inertia acts reducibly): the
  endomorphism $N=0$ and 
\[
\rho_p(f)=\Ind_{W(E_p)}^{W(\QQ_p)}\varkappa,
\]
where $E_p$ is a quadratic extension of $\QQ_p$, and
$\varkappa:W(E_p)^{\ab} \mapsto \CC^\times$ is a quasi-character which
does not factor through the norm map with a quasi-character of
$W(\QQ_p)^{\ab}$ (so that $\rho_p(f)$ is irreducible). Furthermore, if
$\epsilon_p$ denotes the quadratic character corresponding to the
extension $E_p/\QQ_p$, then $\epsilon_p \varkappa = ||\,
. \,||_p^{1-k}$ as quasi-characters of $\QQ_p^\times$.

\item{\bf Supercuspidal Representation II} (inertia acts irreducibly):
  this only happens for $p=2$. $N=0$ and the image of $\rho_p(f)$ is an
  exceptional group.
\end{enumerate}

\begin{remark} In the $n$-dimensional case, the last case occurs only for $p
\le n$.
\label{rmk:1}
\end{remark}

The previous description uses the assumption that the nebentypus of
$f$ is trivial, and the third case relies on the following two facts
due to Henniart (see \cite{Henniart} Theorem 8.2):

\begin{thm}
\label{thm:henniart}
Let $E/F$ be a finite separable extension of degree $n$ of local
fields, and $\rho$ be a linear degree $n$ representations of
$W(E)$. If $R$ denotes its induction to $W(F)$, then:
\begin{enumerate}
\item Let $\varepsilon$ be the character of $F^\times$ that corresponds
  to the determinant of the permutation representation of $W(F)$ acting
  on $W(F)/W(E)$. Then for $x \in F^\times$ one has
\[
\det(R)(x) = \varepsilon(x)^n \det(\rho)(x).
\]
\item Assume that $\rho$ is semi-simple, then $R$ is semi-simple as
  well and one has
\[
\cond(R) = f(E/F) (n\cdot d(E/F) + \cond(\rho)),
\]
where $\cond()$ denotes the exponent of the Artin conductor of the
representation, $f(E/F)$ denotes the inertial degree and $d(E/F)$ is
the discriminant.
\end{enumerate}
\end{thm}

\section{The case of an odd prime}

The previous description and the fact that the unique order two
character of $\ZZ/p^r$ has conductor $p$, gives the well know
condition on the exponents of the distinct types:

\begin{coro} If $f \in S_k(\Gamma_0(N))$ and $p\neq 2$, we have:
\label{coro:valuations}
\begin{enumerate}
\item If $\rho_p(f)$ is principal series then $v_p(N)$ is even.
\item If $\rho_p(f)$ is Steinberg then $\val_p(N)=1$ or $2$.
\item If $\rho_p(f)$ is supercuspidal then $\val_p(N) \ge
  2$. 
\end{enumerate}
Furthermore, in the last case if $E_p/\QQ_p$ is unramified then
$\val_p(N)$ is even.  If $E_p/\QQ_p$ is ramified, then $\val_p(N)$ is
odd unless $\varkappa$ has conductor $1$, $p\equiv 3 \pmod 4$ and
$\varkappa|_{\ZZ_p^\times} = \epsilon_p$, in which case the conductor
is $2$. 
\label{coro:level}
\end{coro}

\begin{proof}
The first two statements are clear. For the last statement, we know
that our representation is the induced representation from a quadratic
extension $E_p$ of $\Q_p$ of a character $\varkappa$. Then
Theorem~\ref{thm:henniart} implies that since $f$ has trivial
nebentypus, we must have
\begin{equation}
\varkappa|_{\ZZ_p^\times} \epsilon_p=1,
\label{eq:determinant}
\end{equation}
where $\epsilon_p$ is the quadratic character that corresponds via
class field theory to the quadratic extension $E_p/\Q_p$.  If such
extension is unramified, the conductor of $\varkappa$ must be
non-zero, since otherwise it will factor through the norm map. Hence
the even condition in the exponent comes from the fact that the
inertial degree is $2$ in this case.

In the case where the extension is ramified, $d(E_p/\Q_p)=1$, hence
Theorem~\ref{thm:henniart} implies that the conductor of the
representation equals $1+\cond(\varkappa)$.  

The conductor of $\epsilon_p$ is $1$, so condition
\eqref{eq:determinant} implies that $\cond(\varkappa) = 1$ and
$\varkappa|_{\Z_p^\times} = \epsilon_p$ or its conductor is even,
which implies the statement. Note that in the first case,
$\varkappa|_{\Z_p^\times}$ is quadratic, and it factors through the
norm map if and only if there is a character of order $4$ in
$\Z_p^\times$ with conductor $p$. This is indeed the case if and only
if $p \equiv 1 \pmod 4$, so the representation is irreducible only
when $p \equiv 3 \pmod 4$.
\end{proof}

\begin{remark}
If the form $f$ has CM by the extension $\Q[\sqrt{-p}]$, then
its local component at $p$ corresponds exactly to $E_p/\Q_p$ being
ramified and $\varkappa$ of conductor $1$.
\end{remark}

\bigskip

Let $\chi$ be the quadratic character associated to the quadratic
extension of $\QQ$ ramified only at $p$. By class field theory, it can
be identified with a character of the id\`ele group, i.e. characters
$\{\chi_q\}_q$, with $\chi_q:\Q_q^\times \mapsto \CC^\times$ satisfying the
following conditions:
\begin{itemize}
\item If $q \neq p$, then $\chi_q$ is unramified, and $\chi_q(q) = \kro{q}{p}$.
\item $\chi_p$ is ramified with conductor $p$, and its value in
  $\ZZ_p^\times$ factors through the unique quadratic character of
  $\FF_p^\times$. Furthermore, $\chi_p(p) = 1$.
\end{itemize}


Given a modular form $f$ in $S_k(\Gamma_0(N))$, we want to study how
the local factors of $f$ change while twisting by $\chi$. Denote by
$\varepsilon_q$ the variation of the local factor of $f$ at $q$ while 
twisting by $\chi_q$, where we choose the same additive character and
Haar measure on both computations.

\begin{remark}
In the correspondence between automorphic forms and representations of
the Weil-Deligne group, twisting an automorphic representation by a
quasi-character, has the effect of twisting the Weil-Deligne representation
by the inverse of the quasi-character, but since our character $\chi_p$ is
quadratic, we can avoid this technical detail.
\end{remark}

\begin{thm}
The number $\varepsilon_q$ is given by:
\begin{enumerate}
\item If $q\neq p$, then $\varepsilon_q = \kro{q}{p}^{\val_q(N)}$.
\item If $\rho_p(f)$ is principal series, then 
\[
\varepsilon_p=
\begin{cases}
\kro{-1}{p} & \text{ if } \val_p(N)\neq 0,\\
\kro{-1}{p}p^k & \text{ if } \val_p(N)=0.
\end{cases}
\]
\item If $\rho_p(f)$ is supercuspidal and $E_p/\QQ_p$ is unramified,
  then $\varepsilon_p = -\kro{-1}{p}$.
\item If $\rho_p(f)$ is supercuspidal and $E_p/\QQ_p$ is ramified, then 
\[
\varepsilon_p = \begin{cases}
1 & \text{ if } \val_p(N)=2,\\
1 & \text{ if } E_p = \QQ_p[\sqrt{p^\star}],\\
-1 & \text{ elsewhere}.
\end{cases}
\]
\item If $\rho_p(f)$ is Steinberg with $\val_p(N)=1$, choosing the
  additive character $\psi$ unramified and the Haar measure normalized
  such that $\int_{\ZZ_p}dx=1$, the local sign is given by
\[
\varepsilon(\rho_p(f),\psi,dx) = \frac{-1}{\chi(p)};
\]
while the local sign of $\rho_p(f) \otimes \chi_p$ is given by
\[
\varepsilon(\rho_p(f)\otimes \chi_p,\psi,dx)= \kro{-1}{p}.
\]
\end{enumerate}
\label{thm:sign}
\end{thm}


\begin{remark}
Although in the second case, the local root number is not just a sign,
the power of the prime $p$ appearing comes from the fact that the
level of the form and its twists are different. 
\end{remark}

\begin{remark}
The result for the Steinberg representation and for the Principal
series when $p\nmid N$ are well known, and can be found for example in
\cite{AL} (Lemma 30 and Theorem 6), although the way to prove it
uses the theory of the Atkin-Lehner involutions as global actions,
while the proof we present is just of local nature.
\end{remark}

The proof of the result is quite elementary, and is mainly based in
the properties that the local constant satisfies, as explained in
\cite{Deligne}. One of the main properties that determine the local
constant is the following:
%

\begin{property}
Let $\rho$ be a virtual $0$-dimensional representation of a finite
extension $E_p/\QQ_p$, then
\[
\varepsilon(\Ind_{W(E_p)}^{W(\QQ_p)}\rho, \psi) = \varepsilon(\rho,\psi \circ \Tr_{E_p/\QQ_p}).
\]
\label{prop:induced}
\end{property}

See \cite{Deligne}, Theorem $4.1$ for a proof of the existence of
local constants with the appropriate $4$ conditions.

\begin{proof}[Proof of Main Theorem]
We consider each case separately:

\medskip

\noindent (1) If $q \neq p$, the character $\chi_q$ is unramified,
  hence by (5.5.1) of \cite{Deligne}, 
\[
\varepsilon(\rho_q(f) \otimes \chi_q,\psi,dx) = \chi_q\left(q^{\val_q(N)+q \dim(\rho_q(f))}\right) \varepsilon(\rho_q(f),\psi,dx).
\]
Since $\chi_q(q) = \kro{q}{p}$ and $\dim(\rho_q(f)) =2$, the statement
follows.\\

\noindent (2) Since $\rho_p(f)$ is reducible in the principal series case, we
  need to see how the local constant of a quasi-character changes under
  twisting by $\chi_p$. Let $a=\cond(\chi_1)$ be the conductor of the
  quasi-character $\chi_1$; chose $\psi$ to be an additive character of
  $\QQ_p$ with $\cond(\psi)=0$ (i.e. $\psi|_{\ZZ_p}=1$ but
  $\psi|_{\frac{1}{p}\ZZ_p}\neq 1$) and the Haar measure $dx$ such
  that $\int_{\ZZ_p} dx =1$. The local epsilon factor is then given by

\[
\int_{\ZZ_p^\times} \chi_1^{-1}\left(\frac{x}{p^a}\right)
\psi\left(\frac{x}{p^a}\right)d\frac{x}{p^a}= \chi_1(p)^a p^a \sum_{b
  \in \Z_p^\times/U_a} \chi_1^{-1}(b)\psi\left(\frac{b}{p^a}\right) \int_{U_a} dx.
\]
The normalization $\int_{\ZZ_p} dx = 1$ implies that
$\int_{\ZZ_p^\times} dx = \frac{p-1}{p}$ and $\int_{U_a}dx =
\frac{1}{p^a}$. Since the conductor of $\psi$ is $0$,
$\psi\left(\frac{1}{p^a}\right) = \exp(2\pi i/p^a)^c$ for some $c$
prime to $p$. Then
\[
G(\chi_1^{-1},c) = \chi_1^{-1}(c) \left(\sum_{b \in \Z_p^\times/U_a}
\chi_1^{-1}(b)\exp\left(\frac{2\pi ib}{p^a}\right)\right),
\]
is a Gauss sum. If we compute the product of the epsilon factor
corresponding to $\chi_1$ and the one corresponding to $\chi_1^{-1}||\,.\,||_p^{1-k}$,
we get that the local factor is given by

\[
||p^{a}||_p^{1-k} G(\chi_1^{-1},c) G(\chi_1,c) =\frac{p^{ak}}{p^a} p^a
\chi_1(-1)=p^{ak} \chi_1(-1).
\]
The middle equality is a classical result of Gauss sums, see for
example \cite{Davenport} (Ex. 13(9), p. 295). If we replace $\chi_1$
by $\chi_1 \chi_p$ in the previous computation, we get that the two
numbers differ by $\chi_p(-1) = \kro{-1}{p}$ and a power of $p$ if the
level of $f$ and that of $f \otimes \chi$ are not equal, as claimed.

\noindent (3) Since $(\Ind_{W(E_p)}^{W(\QQ_p)} \varkappa)\chi_p =
\Ind_{W(E_p)}^{W(\QQ_p)}({\chi}_p \varkappa) $ (where in the second term
of the equality we are considering the restriction of ${\chi}_p$ to
$W(E_p)$), we can  apply the Property stated before to $\rho =
{\chi}_p \varkappa - \varkappa$. Then
\[
\frac{\varepsilon(\Ind_{W(E_p)}^{W(\QQ_p)}{\chi}_p \varkappa,
  \psi)}{\varepsilon(\Ind_{W(E_p)}^{W(\QQ_p)}\varkappa,\psi)} =
\frac{\varepsilon({\chi}_p\varkappa, \psi \circ
  \Tr_{E_p/\QQ_p})}{\varepsilon(\varkappa,\psi \circ \Tr_{E_p/\QQ_p})}.
\]
This allows to restrict to the $1$-dimensional case. Recall that since
$E_p/\Q_p$ is unramified, $\varkappa$ is ramified (as was pointed out
in the proof of Corollary~\ref{coro:level}). Let $\Om_p$ denote the
ring of integers of $E_p$, take $p$ as a local uniformizer, and let
$a = \cond(\varkappa)$ be the conductor of $\varkappa$. Recall that since
the nebentypus of $f$ is trivial, $\varkappa|_{\ZZ_p^\times} = 1$. Let $\psi$
be an additive character with $\cond(\psi)=0$ as before, and $dx$ a Haar
measure such that $\int_{\Om_p}dx=1$. To easy notation we will denote
by $\tilde{\psi}$ the additive character $\psi \circ
\Tr_{E_p/\QQ_p}$. The local factor is given by
\begin{multline*}
\varepsilon(\varkappa,\tilde{\psi},dx)=\int_{\Om_p^\times}
\varkappa^{-1}\left(\frac{x}{p^a}\right)\tilde{\psi}\left(\frac{x}{p^a}\right)d
\frac{x}{p^a} =\\ \varkappa(p)^a p^{2a} \left(\sum_{b \in \Om_p^\times/U_a}
\varkappa^{-1}(b) \tilde{\psi}\left(\frac{b}{p^a}\right)\right) \int_{U_a}dx.
\end{multline*}
The middle sum can be written as
\begin{multline*}
\sum_{\alpha \in (\Om_p/p^a)^\times / (\ZZ_p/p^a)^\times} \sum_{\beta \in (\ZZ_p/p^a)^\times} \varkappa^{-1}(\alpha \beta) \psi\left(\frac{\Tr(\alpha)\beta}{p^a}\right) =\\ \sum_\alpha \varkappa^{-1}(\alpha)\sum_\beta \psi\left(\frac{\Tr(\alpha)\beta}{p^a}\right).
\end{multline*}

If $p^{a-1}\nmid \Tr(\alpha)$, the last sum is zero, since it is a
sum over all primitive roots of unity of order at least $p^2$. For
elements where $p^{a-1}|| \Tr(\alpha)$ (i.e. $p^{a-1}\mid \Tr(\alpha)$
but $p^a \nmid \Tr(\alpha)$), the last sum is $-1$. Such elements are
of the form
\[
r (p^{a-1}+\beta \sqrt{\delta}), \quad r \in (\ZZ_p/p^a)^\times \text{ and }p \nmid \beta.
\]
Modulo multiplication by elements of $(\ZZ_p/p^a)^\times$ in
$\Om_p/p^a$, for $p \nmid s$, we have that
\[
(p^{a-1}+\beta \sqrt{\delta})^s \sim (p^{a-1}+s^{-1} \beta \sqrt{\delta}),
\]
so the sum with these terms is
\[
(-1) \cdot \sum_{p^{a-1}||\Tr(\alpha)} \varkappa^{-1}(\alpha) = (-1) \cdot \sum_{s \in (\ZZ_p/p^a)^\times} \varkappa^{-1}(p^{a-1}+\sqrt{\delta})^s.
\]
Since the conductor of $\varkappa$ is $p^a$, $\varkappa$ is
non-trivial on such elements and the last sum is zero. Then the only
remaining terms are the ones with $\Tr(\alpha)=0$ (i.e. $\alpha =
\sqrt{\delta}$ for a non-square element $\delta$), and in this case
all terms of the last sum are $1$ so we get that

\[
\varepsilon(\varkappa,\tilde{\psi},dx) = p^{a-1}(p-1) \varkappa(\sqrt{\delta})^{-1}\varkappa(p)^a,
\]
where $\delta$ is a non-square in $\QQ_p$. If we make the same
computation with ${\chi}_p \varkappa$, we get that
\[
\varepsilon_p =
\frac{p^{a-1}(p-1)\varkappa(\sqrt{\delta})^{-1}\chi_p(-\delta)\varkappa(p)^a}{p^{a-1}(p-1)\varkappa(\sqrt{\delta})^{-1}\varkappa(p)^a}
= -\kro{-1}{p}.
\]

\medskip

\noindent(4) This case is similar to the previous one, the main
difference is that 
using the commutative diagram
\eqref{eq:commutative1}, we need to compose our character with the
norm map. This gives another character (that abusing notation we also
denote $\chi_p$) which 
satisfies ${\chi}_p|_{\Om_p^\times}=1$, because the conductor of
$\chi_p$ is $p$ and the norm map from $\Om_p^\times$ to $\ZZ_p^\times$
gives only squares modulo $p$. Then the terms in the sum are the same
for $\varkappa$ and ${\chi}_p\varkappa$. Since $E_p/\QQ_p$ is
ramified, the conductor of $\psi \circ \Tr = 1$, hence if we take $\pi
= \sqrt{p\delta}$ as a local uniformizer, the local factor is given by
\[
\varepsilon(\varkappa,\tilde{\psi},dx)=
\varkappa(\pi^{\cond(\varkappa)+1})\int_{\Om_p^\times}
\varkappa^{-1}(x)\tilde{\psi}\left(\frac{x}{\pi^{2r+1}}\right)d
\frac{x}{\pi^{2r+1}}.
\]
The local factor of the twisted representation is given by
\[
\varepsilon(\tilde{\chi}_p\varkappa,\tilde{\psi},dx)=
\chi_p(\Norm(\pi))^{\cond(\varkappa)+1}
\varepsilon(\varkappa,\tilde{\psi},dx).
\]
Hence the quotient equals $1$ if $\cond(\varkappa)=1$ or 
\[
\varepsilon_p=\kro{\Norm(\pi)/p}{p} = \kro{-\delta}{p}.
\]
In particular, if $\delta=\kro{-1}{p}$ up to squares,
$\varepsilon_p=1$, while if it does not, then $\varepsilon_p=-1$ as
claimed.

\medskip

\noindent(5) Comes from the definition of the local epsilon factor
attached to the representations where the nilpotent endomorphism is
not trivial and how the local epsilon factors changes for the
principal series.
\end{proof}

Let $f$ be in $S_k(\Gamma_0(N))$, where  $N = p^r N'$, with $p \nmid
N$, and let $\varepsilon(f)$ be its functional equation sign. Let
$\chi_p$ be as before, and let $f \otimes \chi_p$ denote the newform
obtained while twisting $f$ by $\chi_p$. Denote by $N(f\otimes
\chi_p)$ its level. The previous statement allows the following
classification:
\begin{coro} With the previous notation, we have the following
  computational criteria to compute the local type:
\begin{itemize}
\item $\pi_p(f)$ is \emph{Steinberg} if $\val_p(N)=1$ or $\val_p(N(f \otimes\chi_p))=1$.
\item $\pi_p(f)$ is \emph{Principal Series} if it is not Steinberg, $2
  \mid \val_p(N)$ and
\[
\varepsilon(f \otimes \chi_p) = \chi_p(N') \varepsilon(f)
  \kro{-1}{p}.
\]
\item $\pi_p(f)$ is \emph{Supercuspidal} if it is not of the above
  type. Furthermore, if $\val_p(N)$ is even and greater than $2$,
  $\pi_p(f)$ is induced from the unramified quadratic extension of
  $\Q_p$, while if $\val_p(N)$ is odd and greater than $2$,
\begin{itemize}
\item $\pi_p(f)$ is induced from the extension $\QQ_p[\sqrt{p^\star}]$
if $\varepsilon(f \otimes \chi_p) = \chi_p(N') \varepsilon(f)$.
\item $\pi_p(f)$ is induced from the extension
  $\QQ_p[\delta\sqrt{p^\star}]$ (for any non-square $\delta$) if
  $\varepsilon(f \otimes \chi_p) = -\chi_p(N') \varepsilon(f)$.
\end{itemize}

\end{itemize}
\label{coro:sign}
%

\end{coro}
\begin{remark}
In the case $\val_p(N)=2$, twisting only allows to determine the type, but does
not distinguishes from which quadratic extension the representations
is induced from.
\end{remark}

%




\begin{remark}
We can replace the global functional equation sign in the last two
corollaries by the local Atkin-Lehner involution $W_p$. Then the same
statements are true replacing $\varepsilon()$ by the eigenvalue of
$W_p$ and removing the factor $\chi_p(N')$.
\end{remark}

\begin{remark}
If $f \in S_k(\Gamma_1(N),\epsilon)$ and for $p \mid N$ the character
$\epsilon_p=1$, the same result holds. 
\end{remark}


\section{The case $p=2$}

When $p=2$, there are more representations of the Weil group. It is a
classical result that all subgroups of $\PGL_2(\CC)$ are isomorphic
to: a cyclic group, a dihedral group, $A_4$, $S_4$ or $A_5$. The $A_5$ case
cannot happen since the Galois group $\Gal(\bar{\Q}_p/\Q_p)$ is
solvable. 

The case $p \neq 2$ did not include exceptional cases since for odd
characteristic all $2$-dimensional representations of $W(\QQ_p)$ are
of dihedral type as mentioned in Remark~\ref{rmk:1}. But for $p=2$,
the $A_4$ and $S_4$ case might occur. Weil proved in \cite{Weil} that
over $\Q$, the $A_4$ case actually does not occur, so the only
exceptional case has image $S_4$. Furthermore, there are only $8$
cases with projective image $S_4$ and all cases it corresponds to the field
extension of $\Q_2$ obtained by adding the coordinates of the
$3$-torsion points of the elliptic curves (see also \cite{BR}, Section 8):
\begin{equation}
E_1^{(r)}: ry^2 = x^3+3x+2, \qquad \: \qquad r \in \{\pm 1, \pm 2\},
\end{equation}
and
\begin{equation}
E_2^{(r)}: ry^2 = x^3-3x+1, \qquad \: \qquad r \in \{\pm 1, \pm 2\}.
\end{equation}

Note that there are $3$ quadratic extensions of $\Q$ which ramify only
at $2$. We will denote by $\chi_{-1}, \chi_{2}$ and $\chi_{-2}$ the
quadratic character that corresponds to the such quadratic extensions,
where $\chi_i$ corresponds to $\Q[\sqrt{i}]$ (of conductor $4$ the
first one and $8$ the last two ones). Then the 4 curves of each type
are twists of each other, where $E_i^{(r)} \otimes \chi_j =
E_i^{(rj)}$ (abusing notation and considering the supra-indices modulo
the equivalence relation given by squares). Furthermore, by \cite{Rio}
(Section 6), the level of the modular form is $2^7$ in the case of the
curve $E_1^{(r)}$ (with $r \in \{\pm1, \pm2\}$), $2^4$ for the curve
$E_2^{(1)}$, $2^3$ for the curve $E_2^{(-1)}$ and $2^6$ for the curve
$E_2^{(\pm2)}$.

Before stating the equivalent of Corollary~\ref{coro:valuations},
recall that there are $7$ extensions of $\Q_2$. One of them is
unramified, two of them have discriminant with valuation $2$
(corresponding to $\sqrt{3}$ and $\sqrt{7}$) and four of them have
discriminant with valuation $3$ (corresponding to $\sqrt{2}$,
$\sqrt{10}$, $\sqrt{-2}$ and $\sqrt{-10}$).  With the previous
notations, we have
\begin{coro}
For $p=2$, we have:
\begin{itemize}
\item If $\rho_2(f)$ is principal series then $v_2(N)$ is even (but
  not $2$).
\item If $\rho_2(f)$ is Steinberg then $\val_2(N) \in \{1,4,6\}$.
\item If $\rho_2(f)$ is supercuspidal then $\val_2(N) \ge
  2$. Furthermore, depending on the different extensions we have:
\begin{itemize}
\item If $E_2/\Q_2$ is unramified then $\val_2(N)$ is even and greater
  or equal to $2$.
\item If $E_2/\Q_2$ is ramified with valuation $2$ then $\val_2(N) =5$
  or it is even and greater or equal to $6$.
\item If $E_2/\Q_2$ is ramified with valuation $3$, then $\val_2(N) = 8$
  or it is odd and greater or equal to $9$.
\end{itemize}
\item If $\rho_2(f)$ is supercuspidal of type II, then $\val_2(N) \in
  \{3,4,6,7\}$. 
\end{itemize}
\end{coro}


\begin{proof}
The proof is the same as before for the first two cases. The
supercuspidal dihedral case is also the same for the unramified extension,
while for the ramified extension, note that if
$\varkappa|_{\Z_2^\times} = \epsilon_2$, then the conductor of
$\varkappa$ is $2 \cdot \cond(\epsilon_2) -1$ or it is even (and
greater). Then Theorem~\ref{thm:henniart} implies that the conductor
of the induced representation is $3 \cdot \cond(\epsilon_2)-1$ or
congruent to $\cond(\epsilon_2)$ modulo $2$.
\end{proof}

\begin{remark}
Note that in the supercuspidal dihedral case induced from ramified
quadratic extensions, all the representations are irreducible, since
there are no characters of order $4$ and conductor $4$ nor characters
of order $4$ and conductor $8$.
\end{remark}

\begin{remark}
In the principal series case, the levels $1, 4$ and $6$ are twists of
each other, since all characters of conductor $1, 2$ or $3$ are at
most quadratic.
\end{remark}

We summarize the previous corollary in Table~\ref{table:p2} (where we
used the notation: PS meaning principal series, ST meaning Steinberg,
SCI$a$, SCI$b$, SCI$c$ meaning Supercuspidal of dihedral type and the
index $a,b$ of $c$ meaning induced from an extension with discriminant
$0, 2$ and $3$ respectively; and SCII meaning supercuspidal of second type).

\begin{table}[h]
\begin{tabular}{|c|c|c|c|}
\hline
$\val_2(N)$ & Types & $\val_2(N)$ & Types\\
\hline\hline
$0$ & PS & $6$ & PS, ST, SCI$a$, SCI$b$, SCII \\
\hline
$1$ & ST & $7$ & SCII\\
\hline
$2$ & SCI$a$ & $8$ & PS, SCI$a$, SCI$b$, SCI$c$\\
\hline
$3$ &SCII &  odd $\ge 9$ & SCI$c$\\
\hline 
$4$ &PS, ST, SCI$a$, SCII &even $\ge 10$ & PS, SCI$a$, SCI$b$ \\
\hline 
$5$ & SCI$b$ & & \\
\hline
\end{tabular}
\caption{Possible types for $p=2$. \label{table:p2}}
\end{table}

Let $\chi$ denote the character associated by class field theory to
any of the characters $\chi_i$, $i \in \{-1,\pm 2\}$ and denote by
$\varepsilon_q$ the change of the variation of the local factor at $q$
under twisting by $\chi$. Then $\varepsilon_q$ it is given by:

\begin{thm}
The number $\varepsilon_q$ (up to a power of $2$) is given by:
\begin{enumerate}
\item If $q\neq 2$, then $\varepsilon_q = \chi(q)^{\val_q(N)}$.
\item If $\rho_2(f)$ is principal series, then 
$\varepsilon_2= \chi(-1)$.
\item If $\rho_2(f)$ is supercuspidal and $E_2/\QQ_2$ is unramified,
  then $\varepsilon_2 = -\chi(-1)$.
\item If $\rho_2(f)$ is supercuspidal and $E_2/\QQ_2$ is ramified, then 
\[
\varepsilon_2 = \begin{cases}
1 & \text{ if } \cond(\varkappa) = 2\cdot \cond(\epsilon_2)-1,\\
1 & \text{ if } E_2 \text{ corresponds to }\chi_2,\\
-1 & \text{ elsewhere}.
\end{cases}
\]
\item If $\rho_2(f)$ is Steinberg with $\val_2(N)=1$, choosing the
  additive character $\psi$ unramified and the Haar measure normalized
  such that $\int_{\ZZ_2}dx=1$, the local sign is given by
\[
\varepsilon(\rho_2(f),\psi,dx) = \frac{-1}{\chi(2)};
\]
while the local sign of $\rho_2(f) \otimes \chi_2$ is given by
\[
\varepsilon(\rho_2(f)\otimes \chi_2,\psi,dx)= \chi_2(-1).
\]
\end{enumerate}
\label{thm:signp2}
\end{thm}
\begin{proof}
The proof is almost the same as the odd case, and can easily be checked. 
\end{proof}
\begin{remark}
The computation in the non-dihedral supercuspidal case is straight
forward, since one can compute for each elliptic curve the local root
number and use the relations between the curves under twisting. In
Table  \ref{table:ss} we list the local root numbers of each curve at $2$.

\begin{table}[h]
\begin{tabular}{|c|r|c|r|}
\hline
Curve & Root Number & Curve & Root Number\\
\hline\hline
$E_1^{(1)}$ & $1$ & $E_2^{(1)}$ & $1$\\
\hline
$E_1^{(-1)}$ & $1$ & $E_2^{(-1)}$ & $1$\\
\hline
$E_1^{(2)}$ & $-1$ & $E_2^{(2)}$ & $-1$\\
\hline
$E_1^{(-2)}$ & $1$ & $E_2^{(-2)}$ & $1$\\
\hline
\end{tabular}
\caption{Root numbers in the non-dihedral supercuspidal case. }
\label{table:ss}
\end{table}
\end{remark}

\begin{remark} 
Contrary to the odd case, where the variation of the sign under
twisting allows to compute the exact type of a representation, for
$p=2$ this is no longer the case. The only cases that cannot being
distinguished are those of a Principal Series representation and a
dihedral supercuspidal representation induced from the quadratic
extension $\Q_2(\sqrt{2})$, in the case they are not twists of lower
level so in particular $\val_2(N)$ is even and greater than $8$. This
is so by Theorem~\ref{thm:signp2} since, following the previous
notation, $\chi_{-1}(-1) = -1$, $\chi_{-2}(-1) = -1$ and $\chi_2(-1) =
1$.
\label{rmk:examples}
\end{remark}

\begin{example}
\label{ex:1}
Consider the curve $E768b$ in Cremona's notation. Its quadratic twists
by $\chi_{-1}$, $\chi_{2}$ and $\chi_{-2}$ are the curves $E768h$,
$E768d$ and $E768f$ respectively (an online table of the curves and
their first twists can be found in \cite{tornaria}). Looking at the
local root number at $2$ in such tables, we see that they change by
$-1$, $1$ and $-1$ respectively so we are in the condition of the last
Remark. To see whether we are in the Principal Series case or in the
Supercuspidal one, we can search for the curve in the quaternion
algebra ramified at $2$ and infinity. This can be done by choosing the
correct order in such algebra (see \cite{HPS}) and constructing
ideal representatives for it in order to compute the Brandt
matrices (see \cite{Nico} for an effective  way to construct the ideals). It turns
out that all four curves appear in such algebra (although it is clear
that if one does the others do as well), hence the component at $2$ of
all of them is supercuspidal.
\end{example}
\begin{example}
\label{ex:2}
Consider the elliptic curve $E3840c$ in Cremona's notation. Its
quadratic twists by $\chi_{-1}$, $\chi_{2}$ and $\chi_{-2}$ are the
curves $E3840w$, $E3840n$ and $E3840t$ respectively. Their local root
numbers at $2$ show that we are again in the condition of
Remark~\ref{rmk:examples}. However, this curve does not show up in the
quaternion algebra ramified at $2$, hence it is Principal Series at $2$.
\end{example}

The last two examples show that both cases actually do occur, as was
expected, and in particular proves that by only considering the
variation of the local root number under twisting is not enough to
determine the local factor at the prime $2$.

\section{Some remarks on Hilbert modular forms}

Although in all the previous sections we worked only with classical
modular forms, the correspondence between Weil-Deligne representations
and Hilbert modular forms works just as well. The properties/existence
of the local root numbers do also, so we could just started with a
Hilbert modular form over a totally positive number field $K$ in the
first case. All the local computations are the same but the
problem is that there might be no global Hecke character $\chi$
to twist by. This comes from the fact that a totally positive number
field $K$ (other than $\Q$) does have totally positive units different
from $1$, which does not happen over $\Q$, so for a character
$\chi_{\id{p}}$ to be well defined, it needs to be trivial at totally
positive units.

To overcome this problem, starting from the character $\chi_{\id{p}}$,
we chose an auxiliary prime $\id{q}$ which does not divide the level
of the Hilbert modular form, and such that $\chi_{\id{p}}
\chi_{\id{q}}$ is trivial on totally positive units. Such primes
always exist, since if $\chi_{\id{p}}$ is non-trivial on units, we can
chose a basis $\{\nu_1, \ldots, \nu_r\}$ of the totally positive units
such that $\chi_{\id{p}}(\nu_1)=-1$ and $\chi_{\id{p}}(\nu_i) = 1$ for
all $2 \le i \le r$. This is equivalent to say that our prime $\id{p}$
is inert in the (ring of integers of the) quadratic extension
$K[\sqrt{\nu_1}]$ and splits in the extension $K[\sqrt{\nu_i}]$, for
$2 \le i \le r$. Then any prime $\id{q}$ with the same spliting
behaviour satisfies our hypothesis (and they always exist by
Tchebotarev).

The behaviour of twisting by
$\chi_{\id{q}}$ is computed using the first case of
Theorem~\ref{thm:sign}, where we need to replace the quadratic symbol
by $\chi_{\id{q}}(\pi)$ for $\pi$ a local uniformizer of
$K_{\id{p}}$. In this way, we can extract the information needed to
compute the local factor at $\id{p}$.

\begin{example} Let $K=\Q[\sqrt{5}]$. The group of totally positive
  units is generated by the element $\langle
  \frac{3+\sqrt{5}}{2}\rangle$. Let $\id{P}_{31} =(6+\sqrt{5})$ be a
  prime ideal of norm $31$ in $K$. In this case,
  $\chi_{\id{P}_{31}}\left(\frac{3+\sqrt{5}}{2}\right) = \chi_{31}(14)=1$, so no
  auxiliary prime is needed and everything works as over $\Q$. For
  example, consider the space of weight $(2,2)$ and level
  $\id{P}_{31}^2$ Hilbert modular forms. This can be computed using
  Demb\'el\'e algorithm (see \cite{lassina}) which is implemented in
  Magma. It turns out that there are $3$ forms having $\Q$ as
  coefficient field. One of them is the twist of the elliptic curve of
  conductor $\id{P}_{31}$ given in \cite{lassina2}, Example 1, hence
  both curves are Steinberg at the prime $\id{P}_{31}$.

The other two curves are one twist of the other, and have Weierstrass
equation:
\begin{equation*}
E:y^2+y = x^3 -x^2 -\left(\frac{7+3\sqrt{5}}{2}\right)x,
\end{equation*}
(this equation was computed by Lassina Demb\'el\'e for us) and its
twist has global minimal model:
\begin{equation*}
E_{\id{P}_{31}}: y^2+\sqrt{5}y =
x^3-\left(\frac{1-\sqrt{5}}{2}\right)x^2 -(639+285\sqrt{5})x -\left(\frac{4733+2113\sqrt{5}}{2}\right).
\end{equation*}
If we compute the sign of their L-series, we see that $E$ has sign
$-1$ while $E_{\id{P}_{31}}$ has sign $+1$ (actually using SAGE, \cite{sage},
one can check that $E$ has rank $1$ while $E_{\id{P}_{31}}$ has
rank $0$). Since $\kro{-1}{31} = -1$, we conclude that $E$ is
principal series at $\id{P}_{31}$.
\end{example}

\begin{example} 
\label{ex:4}
Let $K$ be the totally real field of discriminant $257$ obtained by
adding to $\Q$ a root of the polynomial $t^3+2t^2-3t-1$ (it is not a
Galois extension, since the discriminant is a non-square). The units
in the ring of integer are generated by $\<t,t-1>$, with signatures
$(-1,-1,1)$ both of them. The class number of $K$ is $1$, but the ray
class number is $2$. The totally positive units are generated by
$\<t(t-1),t^2>$. Let $\id{P}_3 = \<t+1>$. It is a non-principal ideal
for the ray class group (the sign of $t+1$ under the three embeddings
are $(+,-,+)$). Also, $\chi_{\id{P}_3}(x(x-1)) = \kro{2}{3}=-1$, so it
does not define a global character. The space of parallel
weight $2$ forms of level $\id{P}_3^2$ has dimension $2$, and it splits
into two eigenforms with rational coefficients (this space was
computed to us by John Voight). One form is a twist of the other one
by the narrow class character. One of the forms correspond to the elliptic curve

\[
E: y^2 + (t^2 + 3t + 3)xy + y = x^3 + (t^2 +t - 1)x^2 + (4t^2 +
19t +4)x + (4t^2 + 9t +2).
\]
One way to prove that the curve is modular for the above modular form,
is to notice that the curve has torsion $\Z/6\Z$, so it is modular and
since there are no other forms in the space, it matches one of the two
forms in our space (this argument and the equation for the elliptic
curve is due to John Voight). Now we search for a prime ideal $\id{P}$
such that it has the same sign in the totally positive units as
$\id{P}_3$. A small search reveals that the prime ideal $\id{P}_7 =
\<2t+1>$ satisfies the required property, since
$\chi_{\id{P}_7}(t(t-1))=\kro{-1}{7}=-1$. So we can compute the twist
of $E$ by the ideal $\id{P}_3 \id{P}_7 =\<2t^2+3t+1>$. It is given by
the equation (in global minimal model)
\begin{multline*}
E_{\id{P}_3\id{P}_7}:y^2 + (t+1)xy = x^3 + (t+1)x^2 +\\ (-863t^2-1791t-442)x + (18919t^2+40953t+10179).
\end{multline*}

Its conductor has valuation $1$ at the prime ideal $\id{P}_3$, hence
both modular forms of level $\id{P}_3^2$ are Steinberg at
$\id{P}_3$. Note that in this case, $\id{p}_3^ 2$ is the smallest
conductor of any twist of the curve and has valuation $2$ at
$\id{p}_3$ (so it is not in any table precomputed before).


\end{example}

\bibliographystyle{alpha} 
\bibliography{computinglocaltype}

\end{document}